\documentclass[conference]{IEEEtran}
\IEEEoverridecommandlockouts
\usepackage{graphicx}
\usepackage{array}
\usepackage{textcomp}
\usepackage{cite}
\usepackage{url}
\usepackage{amsmath}
\usepackage{multirow}
\usepackage{bm}
\usepackage{verbatim}

\ifCLASSINFOpdf
\else
\fi
\begin{document}
%
\title{Optimal Investment on Series FACTS Device Considering Contingencies}

\author{\IEEEauthorblockN{Xiaohu Zhang and Kevin Tomsovic}
\IEEEauthorblockA{Department of Electrical Engineering and Computer Science\\
The University of Tennessee, Knoxville\\
Email: xzhang46@vols.utk.edu, tomsovic@utk.edu}
\and
\IEEEauthorblockN{Aleksandar Dimitrovski}
\IEEEauthorblockA{Power and Energy Systems Group\\
Oak Ridge National Laboratory\\
Email: dimitrovskia@ornl.gov}
\thanks{This work was supported in part by ARPAe (Advanced Research Projects
	Agency Energy), and in part by the Engineering Research Center Program
	of the National Science Foundation and the Department of Energy under
	NSF Award Number EEC-1041877 and the CURENT Industry Partnership Program.}}


%


\maketitle

\begin{abstract}
Series FACTS devices are capable of relieving congestion and reducing generation costs in the power system. This paper proposes a planning model to optimally allocate TCSCs in the transmission network considering the base case and $N-1$ contingencies. We consider a single target year and select three distinct load patterns to accommodate the yearly load profile. A reformulation technique is utilized to linearize the nonlinear power flow constraint introduced by a TCSC. Numerical case studies based on the IEEE 118-bus system demonstrate the high performance of the proposed planning model.  

\end{abstract}


%
\IEEEpeerreviewmaketitle

\section*{Nomenclature}
\subsection*{Indices}
\addcontentsline{toc}{section}{Nomenclature}
\begin{IEEEdescription}[\IEEEusemathlabelsep\IEEEsetlabelwidth{$V_1,V_2,V_3$}]
	\item[$c$] Index of states; $c=0$ indicates the base case; $c>0$ is a contingency state.
	\item[$i, \ j$] Index of buses.
	\item[$k$] Index of transmission elements.
	\item[$m$] Index of loads.
	\item[$n$] Index of generators.
	\item[$t$] Index of load levels.
\end{IEEEdescription}

\subsection*{Variables}
\addcontentsline{toc}{section}{Nomenclature}
\begin{IEEEdescription}[\IEEEusemathlabelsep\IEEEsetlabelwidth{$V_1,V_2,V_3$}]
	\item[$P^g_{nct}$] Active power generation of generator $n$ for state $c$ at load level $t$.
	\item[$P_{kct}$] Active power flow on branch $k$ for state $c$ at load level $t$.
	\item[$\Delta P^{g,up}_{nct}$] Active power generation adjustment up of generator $n$ for state $c$ at load level $t$.
	\item[$\Delta P^{g,dn}_{nct}$] Active power generation adjustment down of generator $n$ for state $c$ at load level $t$.
	\item[$\Delta P^d_{mct}$] Load shedding amount of load $m$ for state $c$ at load level $t$.
	\item[$x^V_{k}$] Reactance of a TCSC at branch $k$.
	\item[$\theta_{kct}$] The angle difference across the branch $k$ for state $c$ at load level $t$.
	\item[$\delta_{k}$] Binary variable associated with placing a TCSC on branch $k$.
\end{IEEEdescription} 

\subsection*{Parameters}
\addcontentsline{toc}{section}{Nomenclature}
\begin{IEEEdescription}[\IEEEusemathlabelsep\IEEEsetlabelwidth{$V_1,V_2,V_3$}]
	\item[$a_n^g$] Cost coefficient for generator $n$.
	\item[$a_n^{g,up}$] Cost coefficient for generator $n$ to increase active power.
	\item[$a_n^{g,dn}$] Cost coefficient for generator $n$ to decrease active power.
	\item[$a_{LS}$] Cost coefficient for the load shedding.
	\item[$N_{kct}$] Binary parameter associated with the status of branch $k$ for state $c$ at load level $t$. 
	\item[$P_{nct}^{g,\min}$] Minimum active power output of generator $n$ for state $c$ at load level $t$.
	\item[$P_{nct}^{g,\max}$] Maximum active power output of generator $n$ for state $c$ at load level $t$.
	\item[$P_{mct}^d$]  Active power consumption of demand $m$ for state $c$ at load level $t$.
	\item[$R_n^{g,up}$] Ramp up limit of generator $n$ .
	\item[$R_n^{g,dn}$] Ramp down limit of generator $n$.
	\item[$S_{kct}^{\max}$] Thermal limit of branch  $k$ for state $c$ at load level $t$.
	\item[$\theta_k^{\max}$] Maximum angle difference across branch $k$: $\pi/3$ radians.
	\item[$\pi_{ct}$] Duration of state $c$ at load level $t$.
\end{IEEEdescription} 

\subsection*{Sets}
\addcontentsline{toc}{section}{Nomenclature}
\begin{IEEEdescription}[\IEEEusemathlabelsep\IEEEsetlabelwidth{$V_1,V_2,V_3$}]
	\item[$\mathcal{D}$] Set of loads.
	\item[$\mathcal{D}_i$] Set of loads located at bus $i$.
	\item[$\Omega_{L}$] Set of existing transmission lines.
	\item[$\Omega_{L}^i$] Set of transmission lines connected to bus $i$.
	\item[$\Omega_{T}$] Set of load levels.
	\item[$\Omega_{c}$] Set of contingency operating states.
	\item[$\Omega_{0}$] Set of base operating states.
	\item[$\Omega_{V}$] Set of candidate transmission lines to install TCSC.
	\item[$\mathcal{B}$] Set of buses. 
	\item[$\mathcal{G}$] Set of on-line generators.
	\item[$\mathcal{G}_{i}$] Set of on-line generators located at bus $i$.
\end{IEEEdescription}

\section{Introduction}
\label{introduction}
\IEEEPARstart{C}{ontrol} of power flow in the interconnected power system has increasingly become a major concern for utilities and system operators as a result of power system restructuring and constrained transmission expansion. Unregulated active and reactive power flows may lead to underutilization of transmission line capacity, high losses, loop flows, decreased stability margins, and so on. In general, there are two options for alleviating these problems. The first option is through power system expansion, i.e. building new power plants and transmission lines to relieve congested areas. The second option involves installing modern power flow control devices. The investment cost issues must be taken into consideration for both options; however, the difficulty in obtaining right of way, political obstacles and long construction times are major hurdles for new transmission lines and upgrades. Given these considerations and improvements in power flow control devices, better utilization of existing power system capacities by installing new equipment is increasingly attractive \cite{mybibb:optfacts,mybibb:optTCSCMINLP}.

The air-core series reactor is one technology for power flow control \cite{mybibb:aircore}. The advantage of the series reactor is low cost and the simple control since it only has two states, i.e., switched in or out. A different load distribution  may require a different reactance and so there are limits to the applications of series reactors. Moreover, frequent electromechanical switching shortens the equipment life and can cause system stress. The appearance of Flexible AC Transmission Systems (FACTS) devices in the last two decades provides new opportunities for controlling power flow and maximizing the utilization of the existing transmission lines \cite{mybibb:FACTS1}, \cite{mybibb:SSSC}. A series FACTS device, such as, a thyristor controlled series compensator (TCSC), is capable of continuously varying the impedance so as to control the power flow. By taking advantage of rapid improvements in power electronics technology, FACTS devices offer excellent control and flexibility. In addition, efforts under the Green Electricity Network Integration (GENI) \cite{mybibb:geni} program initiated by the Advanced Research Projects Agency-Energy (ARPA-E), have led to new FACTS like devices \cite{mybibb:Aleks,mybibb:sheng_zheng} with very low cost, which may be commercially available soon.  

Determining the best locations and settings of FACTS devices in a highly interconnected network is a complex task. Evolutionary computation techniques such as genetic algorithm \cite{mybibb:GAFACTS}, differential evolution \cite{mybibb:differentialevoluation}, particle swarm optimization \cite{mybibb:PSO_FACTS} have all been proposed to find the optimal placement of TCSC to enhance the system loadability or minimize the system active power loss. These techniques have some advantages for dealing with non-differential and non-convex problem, but poor scalability and repeatability prevent them from being widely applied. Sensitivity approaches, e.g.,\cite{mybibb:TCSCsens,mybibb:TCSC_ali}, are another group of methods for optimizing TCSC. In \cite{mybibb:TCSCsens}, the optimal locations of TCSC are computed by using the sensitivity of the transfer capability with respect to the line reactance. In \cite{mybibb:TCSC_ali}, an index called single contingency sensitivity (SCS) is introduced. This weighted index indicates the effectiveness of a given branch in alleviating the congestion under all considered contingencies. However, the sensitivity approach cannot ensure optimality. With the advances in branch-and-bound algorithms, mixed-integer programming (MIP) have been employed as well. The authors in \cite{mybibb:optTCSCMINLP} propose an MINLP model to determine the locations and settings for TCSC to enhance the system loadability in a deregulated power market. In \cite{mybibb:TCSCLFB}, the locations and settings of TCSC are identified via MILP and MIQP using the line-flow-based equation (LFB) proposed by \cite{mybibb:LFB}. To eliminate the quadratic terms in the constraints, one variable in the quadratic term is replaced by its respective limit. The phase angle constraint which is essential in the meshed network is not included in the planning model, which limits the model application. In \cite{mybibb:Tao_Ding,mybibb:asu_FACTS1}, benefits of FACTS devices on the economic dispatch (ED) problem was investigated. The bilinear term of the product between the variable reactance and bus voltage angle was linearized by using the big-M method. The nonlinear programming model was reformulated to an MILP model which can be solved by commercial solvers to achieve the global optimums.    

Few of the previously mentioned work have addressed the economic benefits of FACTS devices considering contingency analysis in the planning model. According to \cite{mybibb:opf_redispatch_facts}, the FACTS devices can be utilized in the security-driven redispatching procedure to reduce the generation rescheduling and load shedding cost in the contingencies. Therefore, considering the cost effects brought by the FACTS devices for the contingencies in the allocation problem should provide a more useful investment plan for the system designers.   

This paper proposes an optimization model to optimally allocate TCSC considering the base case and a series of $N-1$ transmission contingencies. To identify the optimal investment on TCSC, a single target year is considered and three distinct load patterns which represent the peak, normal and low load levels are selected to accommodate the yearly load profile. The contingencies may occur in any of the load levels and their contribution to the total planning cost are weighted by their duration in the target year.       

The remaining sections are organized as follows. In Section \ref{reformulation}, a reformulation technique is proposed to linearize the nonlinear power flow constraints introduced by TCSC. Section \ref{optimization_model} illustrates details of the planning model. In Section \ref{case_studies}, the IEEE 118-bus system is selected for the case studies. Conclusions are given in Section \ref{conclusion}.

\section{Reformulation Technique}
\label{reformulation}
\begin{figure}[!htb]
	\centering
	\includegraphics[width=0.4\textwidth]{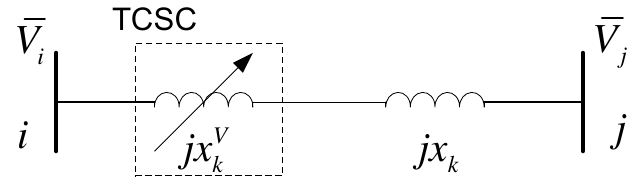}
	\caption{Static representation of TCSC in DCPF.}
	\label{tcsc_steady}
\end{figure}

The steady state model of the TCSC can be represented by a variable reactance $x_k^V$ in series with the transmission line reactance $x_k$, as shown in Fig. \ref{tcsc_steady}.

The total susceptance of the transmission line becomes:
\begin{equation}
b_k^{\prime}=-\frac{1}{x_k+x_k^{V}}=-(b_{k}+b_k^{V})
\end{equation}
where 
\begin{align}
&b_{k}=\frac{1}{x_k}   \\
&b_k^{V}=-\frac{x_k^{V}}{x_k(x_k+x_k^{V})}
\end{align}

In this work, the compensation range of the TCSC is allowed to vary from -70\% to +20\% of its corresponding line reactance\cite{mybibb:TCSC_cost_recovery}. Then the range of the variable susceptance $b^V_{k}$ is:
\begin{align}
b_{k,V}^{\min}&=-\frac{0.2x_k}{x_k(x_k+0.2x_k)}=-\frac{1}{6x_k}  \label{bCVSRmin} \\
b_{k,V}^{\max}&=-\frac{-0.7x_k}{x_k(x_k-0.7x_{k})}=\frac{7}{3x_k}  \label{bCVSRmax}
\end{align}

By introducing a binary variable $\delta_k$ which flags the installation of TCSC on a specific transmission line, the active power flow on that line assuming a DC power flow model can be expressed as:
\begin{align}
&P_k=(b_{k}+\delta_kb_k^{V})\cdot \theta_k  \label{power_CVSR}  \\
&b_{k,V}^{\min} \le b_k^{V} \le b_{k,V}^{\max}   \label{bcvsr_range} 
\end{align}

The nonlinearity in (\ref{power_CVSR}) results from the trilinear term $\delta_kb_k^{V}\theta_k$. To linearize the nonlinear part, a new variable $w_k$ is introduced as:
\begin{equation}
w_k=\delta_kb_k^{V}\theta_k  \label{wij1}
\end{equation}

Then the active power flow equation (\ref{power_CVSR}) becomes:
\begin{equation}
P_k=b_{k}\theta_k+w_k  \label{PCVSR_wij1}
\end{equation}

We multiply each side of the constraint (\ref{bcvsr_range}) with the binary variable $\delta_k$ and combine with $w_k$:
\begin{equation}
\delta_k b_{k,V}^{\min} \le \frac{w_k}{\theta_k}=\delta_kb_k^{V} \le \delta_k b_{k,V}^{\max}  \label{if_ineq}
\end{equation}

The allowable range for $w_k$ is determined by the sign of $\theta_k$:
\begin{equation}
\left\lbrace 
\begin{aligned}
\delta_k\theta_k b_{k,V}^{\min} \le w_k \le \delta_k\theta_k b_{k,V}^{\max}, \ &\text{if} \ \theta_k>0  \\
w_k=0,\ \ \ \ \ \ \ \ \ \ \ \ \ \ \ \ \ \ \ \ \ \ \ \ \ &\text{if} \ \theta_k=0   \\
\delta_k\theta_k b_{k,V}^{\max} \le w_k \le \delta_k\theta_k b_{k,V}^{\min}, \ &\text{if} \ \theta_k<0 
\end{aligned} \right.
\end{equation}

To realize the ``if" constraints, an additional binary variable $y_k$ and the big-M complementary constraints \cite{mybibb:Tao_bigM} are introduced:
\begin{equation}
-M_ky_k+\delta_k\theta_kb_{k,V}^{\min} \le w_k 
\le \delta_k\theta_kb_{k,V}^{\max}+M_ky_k \label{if_1} \\
\end{equation}
\begin{align}
-M_k(1-y_k)+\delta_k\theta_k&b_{k,V}^{\max} \le w_k \nonumber  \\ 
&\le \delta_k\theta_kb_{k,V}^{\min}+M_k(1-y_k) \label{if_2} 
\end{align}

Only one of the two constraints (\ref{if_1}) and (\ref{if_2}) will become active and the other one will be a redundant constraint that is always satisfied with a sufficiently large number $M_k$. Specifically, when $\theta_k< 0$, $y_k$ will be equal to one and the constraint (\ref{if_2}) will be active; when $\theta_k > 0$, $y_k$ will be equal to zero and the constraint (\ref{if_1}) will be active; when $\theta_k=0$, one of these two constraints will drive $w_k$ to be zero regardless of the value of $y_k$. Note that numerical problems may occur if $M_k$ is chosen to be too large \cite{mybibb:TS_contingency}. In this work, $M_k$ is selected as $|\frac{7}{3x_k}\theta_k^{\max}|$.

Now in the constraints (\ref{if_1}) and (\ref{if_2}), the term $\delta_k\theta_k$ is nonlinear which involves the product of a binary variable and a continuous variable. We introduce another new variable $z_k=\delta_k\theta_k$ and linearize using the method in \cite{mybibb:WP_Adams}:
\begin{align}
&-\delta_k\theta_k^{\max} \le z_k \le \delta_k\theta_k^{\max} \label{z1} \\
&\theta_k-(1-\delta_k)\theta_k^{\max} \le z_k \le \theta_k+(1-\delta_k)\theta_k^{\max} \label{z2}
\end{align}

We then substitute $\delta_k\theta_k$ with $z_k$ in the inequalities (\ref{if_1}) and (\ref{if_2}):
\begin{align}
&-M_ky_k+z_kb_{k,V}^{\min} \le w_k 
\le z_kb_{k,V}^{\max}+M_ky_k \label{if3} \\
&-M_k(1-y_k)+z_kb_{k,V}^{\max} \le w_k\le z_kb_{k,V}^{\min}+M_k(1-y_k)\label{if4} 
\end{align}

Thus, the power flow equations (\ref{power_CVSR}) and (\ref{bcvsr_range}) in MINLP format is transformed to an MILP model involving (\ref{PCVSR_wij1}), (\ref{z1})-(\ref{if4}).

\section{Problem Formulation}
\label{optimization_model}
\subsection{Investment Cost for TCSC}
According to \cite{mybibb:tcsc_cost,mybibb:GA_FACT_market,mybibb:PSO_SQP_FACTS}, the investment cost of the TCSC is dependent on its operation range and can be expressed as:
\begin{equation}
I_T=0.0015S_T^2-0.713S_T+153.75     \label{capital_cost}
\end{equation}
where $I_T$ is the cost in $\$/kVar$ and $S_T$ is the installed capacities of TCSC in $Mvar$. We use the capital recovery factor to convert the total investment cost into the annual cost:
\begin{equation}
AI_T=I_T \cdot S_T \cdot 1000 \cdot \frac{d(1+d)^{LT}}{(1+d)^{LT}-1}   \label{annual_cost}
\end{equation}
$d$ is the interest rate and $LT$ is the life time of the FACTS devices. In this work, $d$ is assumed to be $5\%$ and $LT$ is 5 years.

\subsection{Objective Function}
The total planning cost for a single target year includes: 1) operating cost under normal states; 2) operating cost under contingency states; 3) investment cost of the TCSC. It can be formulated as:
\begin{equation}
\min \ \ \sum_{t\in \Omega_T}(\pi_{0t}C_{0t}+\sum_{c\in \Omega_c}\pi_{ct}C_{ct})+\sum_{k\in \Omega_V}AI_T\delta_k  \label{objective}
\end{equation}

In (\ref{objective}), $C_{0t}$ is the operation cost for the normal state under load level $t$, which can be further expressed as:
\begin{equation}
C_{0t}=\sum_{n\in \Omega_{\mathcal{G}}}a_{n}^gP^g_{n0t}  \label{obj_base}
\end{equation}
$a^g_{n}$ is the marginal cost for generator $n$. $C_{ct}$ is the operation cost for the contingency state $c$ under load level $t$, it includes three terms:
\begin{align}
C_{ct}=&\sum_{n\in \Omega_{\mathcal{G}}}a_{n}^gP^g_{nct}+\sum_{m\in \Omega_D}a_{LS}\Delta P^d_{mct}    \nonumber \\
&+\sum_{n\in \Omega_{\mathcal{G}}}(a_n^{g,up}\Delta P^{g,up}_{nct}+a_n^{g,dn}\Delta P^{g,dn}_{nct})                        \label{obj_cont}     
\end{align}
Specifically, the first term is the generator fuel cost under each contingency; the second term is the cost associated with the involuntary load shedding; the third term is the generator rescheduling cost, which indicates that any change from the base operating condition should have a payment to the agents involved.

Note that the total operating hours for a target year is 8760:
\begin{equation}
\sum_{t\in \Omega_T}\pi_{0t}+\sum_{t\in \Omega_T}\sum_{c\in \Omega_c}\pi_{ct}=8760    \label{operating_hour}
\end{equation}

\subsection{Constraints}
With the reformulation, the active power flow through the transmission lines is formulated as:
\begin{align}
& P_{kct}-b_{k}\theta_{kct}+M_{k}'(1-N_{kct}) \ge 0,\ k\in \Omega_{L}\backslash \Omega_{V} \label{norm_c1} \\
& P_{kct}-b_{k}\theta_{kct}-M_{k}'(1-N_{kct}) \le 0,\ k\in \Omega_{L}\backslash \Omega_{V} \label{norm_c2} \\
& P_{kct}-b_{k}\theta_{kct}-w_{kct}+M_{k}'(1-N_{kct}) \ge 0,\ k\in \Omega_{V} \label{CVSR_c1}  \\
& P_{kct}-b_{k}\theta_{kct}-w_{kct}-M_{k}'(1-N_{kct}) \le 0,\ k\in \Omega_{V} \label{CVSR_c2}   \\
&-M_{k}y_{kct}+z_{kct}b_{k,V}^{\min} \le w_{kct}\le z_{kct}b_{k,V}^{\max}+M_{k}y_{kct}, \nonumber \\
&\ \ \ \ \ \ \ \ \ \ \ \ \ \ \ \ \ \ \ \ \ \ \ \ \ \ \ \ \ \ \ \ \ \ \ \ \ \ \ \ \ \ \ \ \ \ \ \ \ \ \ \ k\in \Omega_{V} \label{big_M1}\\
&-M_{k}(1-y_{kct})+z_{kct}b_{k,V}^{\max} \le w_{kct}   \nonumber   \\ 
&\ \ \ \ \ \ \ \ \ \ \ \ \ \ \ \ \ \le z_{kct}b_{k,V}^{\min}+M_{k}(1-y_{kct}),\ k\in \Omega_{V} \label{big_M2} \\
& -N_{kct}\delta_{k}\theta_{k}^{\max} \le z_{kct} \le N_{kct}\delta_{k}\theta_{k}^{\max},\ k \in \Omega_{V} \label{bilinear_c1} \\
& N_{kct}(\theta_{kct}-(1-\delta_{k})\theta_{k}^{\max}) \le z_{kct}   \nonumber \\
&\ \ \ \ \ \ \ \ \ \ \ \ \ \ \le N_{kct}(\theta_{kct}+(1-\delta_{k})\theta_{k}^{\max}),\ k \in \Omega_{V} \label{bilinear_c2}
\end{align}
Constraints (\ref{norm_c1})-(\ref{bilinear_c2}) hold over the sets $\forall c \in \Omega_c\cup \Omega_0, t \in \Omega_T $.

To include the $N-1$ transmission contingency in the optimization model, a binary parameter $N_{kct}$ is introduced to represent the status of the transmission element $k$ in state $c$ under the load condition level $t$. Constraints (\ref{norm_c1}) and (\ref{norm_c2}) are the active power on the lines without TCSC while constraints (\ref{CVSR_c1}) and (\ref{CVSR_c2}) denote the active power on the candidate lines to install TCSC. If $N_{kct}=1$, the DC line flow equations are forced to hold; otherwise, if the line is out of service, a sufficiently large disjunctive factor $M_k'$ ensures that the line flow constraints are not binding. Constraints (\ref{big_M1})-(\ref{bilinear_c2}) are associated with the expansion of the reformulation technique considering multiple operating states and load condition.  

The power balance constraint at each bus is given by:
\begin{align}
&\sum_{n \in \mathcal{G}_i}P^g_{nct}-\sum_{m \in \mathcal{D}_i}(P^d_{mct}-\Delta P^d_{mct})=\sum_{k\in \Omega_{L}^i } P_{kct}  \nonumber \\ 
&\ \ \ \ \ \ \ \ \ \ \ \ \ \ \ \ \ \ \ \ \ \ \ \ \ \ \ \ \ \ \ \ i \in \mathcal{B}, c \in \Omega_c \cup \Omega_0, t \in \Omega_T            \label{P_balance}   \\
&0 \le \Delta P^d_{mct} \le P^d_{mct}, \ m \in \mathcal{D}, c \in \Omega_c, t \in \Omega_T \label{load_shedding}
\end{align}
Constraint (\ref{load_shedding}) ensures that load shedding should not exceed the amount of existing load. 

The system physical limits include the following:
\begin{align}
&\ \ \ \ -N_{kct}S_{kct}^{\max} \le P_{kct} \le N_{kct} S_{kct}^{\max}, \label{Slimit_E} \\
&\ \ \ \ P_{nct}^{g,\min}\le P_{nct}^g \le P_{nct}^{g,\max}, \label{Pg_limit1} \\
&\ \ \ \ \theta_{ref}=0     \label{reference_angle}
\end{align}
Constraints (\ref{Slimit_E})-(\ref{reference_angle}) hold over the sets $\forall c \in \Omega_c\cup \Omega_0, t \in \Omega_T, k \in \Omega_L, n \in \mathcal{G}$.

Constraint (\ref{Slimit_E}) enforces that the power flow on the line is zero when the line is out of service; otherwise, the line flow is bounded by its thermal rating. It should be noted that the short term thermal rating is selected for the contingency state, which is 1.1 times of the normal thermal rating. Constraint (\ref{Pg_limit1}) indicates that power generation is limited by the capacity of the generator. Constraint (\ref{reference_angle}) sets the phase angle of the reference bus to be zero. 

During a certain contingency, the generator redispatching should be based on the corresponding base operating condition. In addition, the adjustments must be within the generators ramp limits. These constraints are formulated as:
\begin{align}
&P_{nct}^g=P_{n0t}^g+\Delta P_{nct}^{g,up}-\Delta P_{nct}^{g,dn}    \label{ramp_gen} \\
&0 \le \Delta P_{nct}^{g,up} \le R_n^{g,up}  \label{ramp_up}   \\
&0 \le \Delta P_{nct}^{g,dn}  \le R_n^{g,dn}  \label{ramp_dn}
\end{align}  
Constraints(\ref{ramp_gen})-(\ref{ramp_dn}) hold over the sets $\forall c \in \Omega_c, t \in \Omega_T, n \in \mathcal{G}$.

\begin{table*}[t]
	\centering
	\caption{Operating Cost in Different States for Peak Load Level and Normal Load Level}
	\label{peak_normal}
	\begin{tabular}{|c|c|c|c|c|c|c|c|c|c|c|}
		\hline
		Cont.&\multicolumn{6}{c|}{Peak Load Level}&\multicolumn{4}{c|}{Normal Load Level}  \\
		\cline{2-11}
		Line&\multicolumn{2}{c|}{Generation Cost}&\multicolumn{2}{c|}{Rescheduling}&\multicolumn{2}{c|}{Load Shedding}&\multicolumn{2}{c|}{Generation Cost }&\multicolumn{2}{c|}{Rescheduling}   \\
		$i-j$&\multicolumn{2}{c|}{($\$/h$)}&\multicolumn{2}{c|}{Cost ($\$/h$)}&\multicolumn{2}{c|}{Cost ($\$/h$)}&\multicolumn{2}{c|}{ ($\$/h$)}&\multicolumn{2}{c|}{Cost ($\$/h$)}   \\
		\cline{2-11}
		&w/o TCSC&w/t TCSC&w/o TCSC&w/t TCSC&w/o TCSC&w/t TCSC&w/o TCSC&w/t TCSC&w/o TCSC&w/t TCSC    \\
		\hline
		Base case&167653&157601&-&-&-&-&123706&122368&-&-    \\
		\hline
		60-61&159993&157601&4767&0&0&0&123706&122368&0&0    \\
		\hline
		8-5&175115&180066&19707&21899&827578&296846&125703&124923&2816&6038    \\
		\hline
		45-49&160467&157601&3974&0&0&0&123706&122368&0&0    \\
		\hline
		5-11&160621&157601&3752&0&0&0&123706&122368&0&0    \\
		\hline
		4-11&160636&157601&3731&0&0&0&123706&122368&0&0    \\
		\hline
		15-19&160627&157601&3723&0&0&0&123706&122368&0&0    \\
		\hline
		47-69&160444&157601&4010&0&0&0&123706&122368&0&0    \\
		\hline
		15-17&160301&157601&7132&0&0&0&123706&122368&0&0    \\
		\hline
		30-17&178805&179052&9616&19526&416031&173271&125042&123585&4804&7772    \\
		\hline
		38-37&165276&157106&13088&3910&29591&0&123503&121929&821&3722    \\
		\hline
		47-49&160627&157601&3723&0&0&0&123706&122368&0&0    \\
		\hline
		17-31&162253&157256&2519&2795&0&0&123706&121919&0&1632    \\
		\hline
		48-49&160642&157601&3719&0&0&0&123706&122739&0&2018    \\
		\hline
		26-30&181362&173580&11348&15843&100861&0&125107&124946&4983&8677    \\
		\hline
		25-27&179328&174383&10674&18675&637935&63954&140070&125404&14763&9431    \\
		\hline
	\end{tabular}
\end{table*}
\section{Numerical Case Studies}
\label{case_studies}
The proposed planning model is applied to the IEEE 118-bus system. The system data can be found in MATPOWER software package \cite{mybibb:MATPOW}. Since only one load pattern is given in the MATPOWER package, we treat that load pattern as the normal load level. The peak load level is 20\% higher than the normal load level and the low load level is 20\% lower. All the simulations are performed on a computer with an Intel(R) Xeon(R) CPU @ 2.40 GHz and 48 GB of RAM. The complete MILP problem is modeled with YALMIP toolbox \cite{mybibb:YALMIP} and the CPLEX solver \cite{mybibb:CPLEX} is selected to solve the model.

\subsection{IEEE 118-Bus System}
The IEEE-118 bus system includes 118 buses, 177 transmission lines, 9 transformers and 19 generators. The total load for the normal case is 4108 MW. 
In a real power system, the number of candidate locations to install power flow control devices is not large due to some practical consideration such as the available space for the substation, limited budget, etc. Therefore, a preliminary study based the sensitivity approach presented in \cite{mybibb:TCSCsens} is performed to obtain candidate locations for TCSC.  Here, 30 transmission lines are selected to be the candidate locations based on this approach. We consider 15 contingencies which significantly affect the total planning cost so the number of operating states is 48.

The planning model suggests that 3 TCSCs should be installed on lines (17-31), (26-30) and (22-23) to reduce the total planning cost. In Table \ref{peak_normal}, the hourly operating cost of each state for the peak load level and normal load level is provided. We categorize the operating cost into three groups: 1) generation cost; 2) rescheduling cost; 3) load shedding cost. The load shedding cost for the normal load level is not included in the table because there is no load shedding for all the contingencies under the normal load condition. It can be seen that the generation cost reduction occurs in the majority of the operating states with the installation of TCSCs. For example, during the peak load level, the generation cost for the normal operating state decreases from 167653 $\$/h$ to 157601 $\$/h$. In contingencies (8-5) and (30-17), the generation cost with TCSCs are higher than the cost without TCSCs. However, a significant reduction of load shedding cost can be observed for these two contingencies. During the normal load level, the generation cost decreases for all the operating states with TCSC but the amount of reduction is not as much as that during the peak load condition.       

\begin{figure}[!htb]
	\centering
	\includegraphics[width=0.45\textwidth]{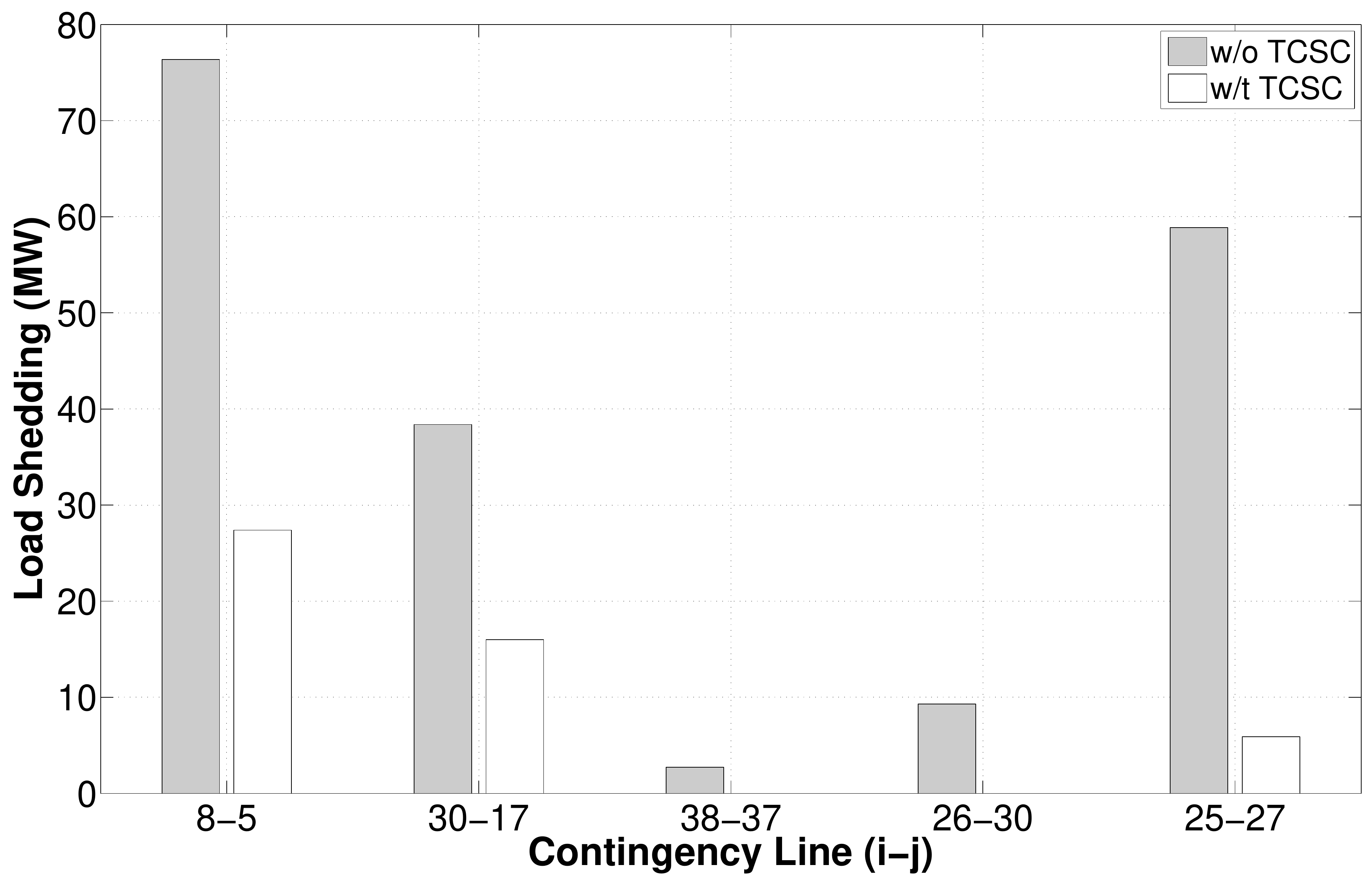}
	\caption{Total load shedding for peak load level in different contingencies.}
	\label{LS}
\end{figure}

Fig. \ref{LS} shows the total load shedding amount for peak load level in contingency (8-5), (30-17), (38-37), (26-30) and (25-27). It can be seen that the load shedding in contingency (38-37) and (26-30) can be eliminated by the TCSCs. For the rest of the three contingencies, the amount of load shedding dramatically decreases.

\begin{figure}[!htb]
	\centering
	\includegraphics[width=0.45\textwidth]{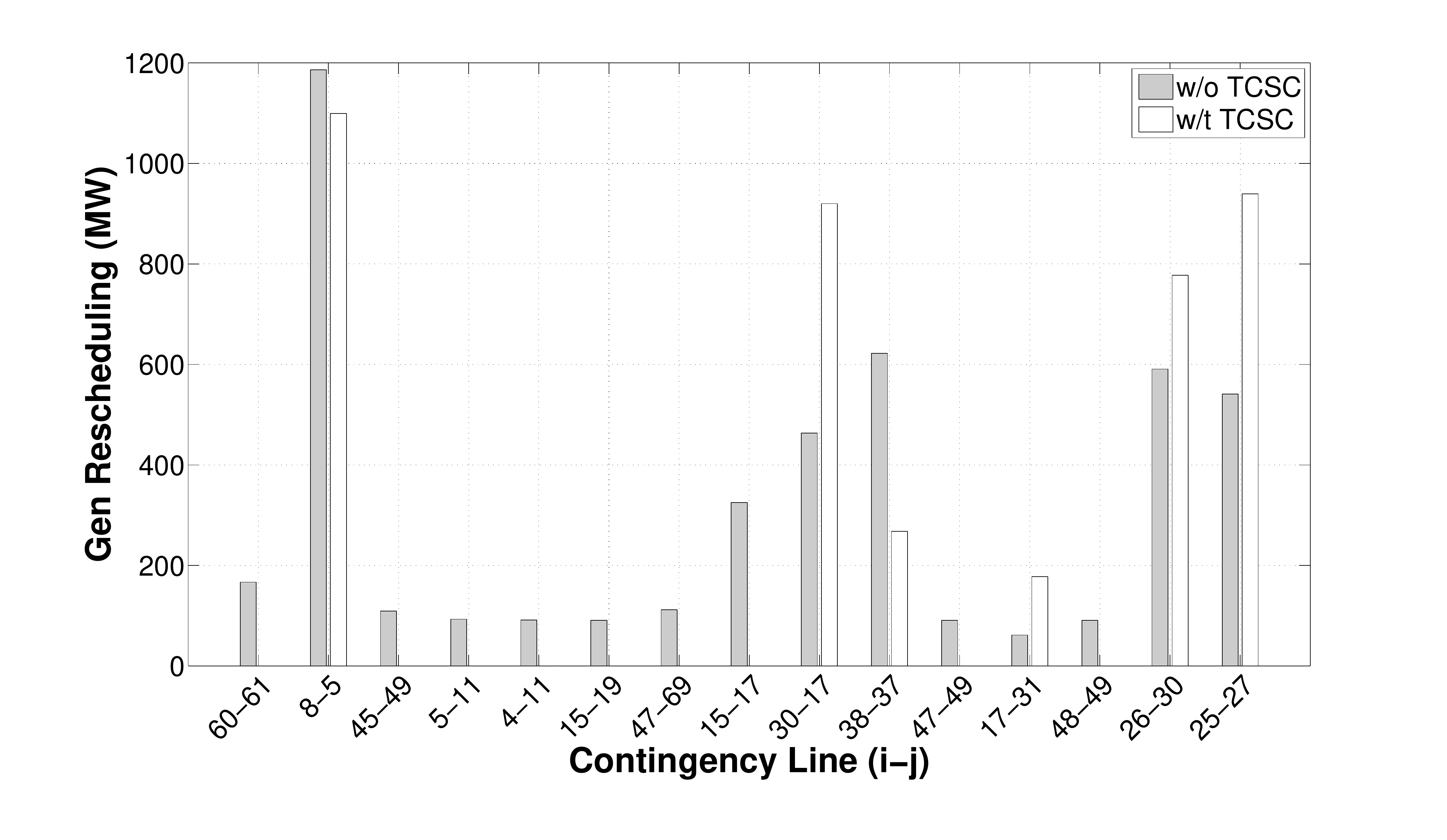}
	\caption{Total generation rescheduling for peak load level in different contingencies.}
	\label{resch}
\end{figure}

The total generation rescheduling amount for each contingency is provided in Fig. \ref{resch}. With three TCSCs, the amount of redispatched power decrease in most of the contingencies but increases for contingencies (30-17), (17-31), (26-30) and (25-27). Three out of these four contingencies involve significant load shedding reduction with TCSC, indicating that TCSC relieves the system congestion and enable more power to be delivered during these contingencies.    
 
\begin{table}[!htb]
	\centering
	\caption{Annual Planning Cost With and Without TCSC}
	\label{annual_saving}
	\begin{tabular}{|c|c|c|}
		\hline
		\multirow{2}{*}{Cost Category}&\multicolumn{2}{c|}{Annual Cost [million \$]}     \\
		\cline{2-3}
		&w/o TCSC&w/t TCSC      \\
		\hline
		Generation cost in normal state&1081.77&1053.48   \\
		\hline
		Generation cost in contingency&33.53&33.07     \\
		\hline
		Rescheduling cost&0.78&0.76           \\
		\hline
		Load shedding cost&8.81&2.34      \\
		\hline
		Investment on TCSC&-&2.31     \\
		\hline
		Total cost&1124.89&1091.96      \\
		\hline
		
	\end{tabular}
\end{table}

Table \ref{annual_saving} provides the annual planning cost for the case with and without TCSC. The installation of three TCSCs decrease the cost in all cost categories. Although the investment on TCSC costs about 2.31 M\$, the total planning cost has an annual reduction of about 33 M\$. This accounts for approximately 2.93\% of the annual planning cost. The computation time required for solving the considered problem is about 1837 s, which is acceptable for the scale of the system considered here.

\section{Conclusion and Future Work}
\label{conclusion}
In this paper, a planning model to optimally allocate TCSC considering base case and a series of $N-1$ contingencies is proposed. The nonlinear power flow constraint introduced by the TCSC is linearized by a reformation technique. Numerical case studies based on IEEE 118-bus system demonstrate the efficiency of the proposed model. Simulation results show that the installation of the TCSC can decrease the generation cost both in the normal operating states and under contingencies. In addition, the generation rescheduling cost and load shedding cost in the contingencies can be reduced. 
Future work is needed to relieve the computational burden and apply the model to a practical large scale power system. Benders Decomposition may serve as one appealing approach to reduce computations.

\bibliographystyle{IEEEtran}
\bibliography{IEEEabrv,mybibb}

\end{document}